\documentclass[10pt]{article}
\usepackage[all]{xy}
\usepackage{amsfonts,amsmath,oldgerm,amssymb,amscd}
\newcommand{\ra}{\rightarrow}

\newcommand{\by}[1]{\stackrel{#1}{\ra}}

\newcommand{\ol}{\overline}          \newcommand{\wt}{\widetilde}
\newcommand{\iso}{\by \sim}

\newtheorem{theorem}{Theorem}[section]
\newtheorem{proposition}[theorem]{Proposition}
\newtheorem{lemma}[theorem]{Lemma}

\newcommand{\BC}{\mbox{$\mathbb C$}}     
     
     \newcommand{\BH}{\mbox{$\mathbb H$}}

     \newcommand{\BR}{\mbox{$\mathbb R$}}

     \newcommand{\BZ}{\mbox{$\mathbb Z$}}

\newcommand{\CC}{\mbox{$\mathcal C$}}     
     
     \newcommand{\CH}{\mbox{$\mathcal H$}}

\newcommand{\op}{\mbox{$\oplus$}}     
\newcommand{\Spec}{\mbox{\rm Spec\,}}
     \newcommand{\Hom}{\mbox{\rm Hom}}

          \newcommand{\ot}{\mbox{$\otimes$}}

\begin{document}

\begin{center}
{\Large \bf $K_0$ of hypersurfaces defined by $x_1^2+\ldots +x_n^2=\pm 1$
}\\
\vspace{.2in} {\large Manoj K. Keshari\footnote{Supported by BOYSCAST
        Fellowship of Department of Science and Technology, India} and
        Satya Mandal$^2$}\\
\vspace{.1in} {\small $^1$Department of Mathematics, IIT Mumbai,
Mumbai - 400076, India \\$^2$Department of Mathematics, University of
Kansas, 1460 Jayhawk Blvd, Lawrence, KS 66045\\
keshari@math.iitb.ac.in, mandal@math.ku.edu}
\end{center}

{\small
\noindent {\bf Abstract:} Let $k$ be a field of characteristic
$\ne 2$ and let
$Q_{n,m}(x_1,\ldots,x_n,y_1,\ldots,y_m)=x_1^2+\ldots+x_n^2-(y_1^2+\ldots+y_m^2)$
be a quadratic form over $k$. Let
$R(Q_{n,m})=R_{n,m}=k[x_1,\ldots,x_n,y_1,\ldots,y_m]/(Q_{n,m}-1)$. In
this note we will calculate $\wt K_0(R_{n,m})$ for every $n,m \geq
0$. We will also calculate $CH_0(R_{n,m})$ and the Euler class group
of $R_{n,m}$ when $k=\BR$.}

\section{Introduction}

Let $A_{n,\BR}=\BR[x_1,\ldots,x_n]/(\sum_1^n x_i^2-1)$ be the
coordinate ring of the real sphere $S^{n-1}$. Then it is well known
(see \cite{ABS}) that $\wt K_0(A_{n,\BR})$ is periodic of period
$8$. More precisely, $\wt K_0(A_{n,\BR})$ is $\BZ,\BZ/2\BZ$ or $0$
depending on whether $n$ is $\{1,5\}$ modulo $8$, $\{2,3\}$ modulo $8$
or $\{0,4,6,7\}$ modulo $8$.  Similarly if
$B_n=\BC[x_1,\ldots,x_n]/(\sum_1^n x_i^2-1)$ is the coordinate ring of
complex $n-1$ sphere, then $\wt K_0(B_n)$ is periodic of period
$2$. More precisely, $\wt K_0(B_n)$ is $\BZ$ or $0$ depending on
whether $n$ is odd or even.

It will be interesting to know if $\wt K_0(A_{n,k})$ is also periodic
for arbitrary field $k$. Further if $\wt
A_{n,\BR}=\BR[x_1,\ldots,x_n]/(\sum_1^n x_i^2+1)$, then we would like
to know if $\wt K_0(\wt A_{n,\BR})$ is periodic. In this paper we answer
these questions in affirmative. 

Some experts may consider these results as easy computations.
However, there is no written reference to these results.  These
results are derived by application of the celebrated results of Swan
\cite{swan1}. We are confident that this article will serve as valuable
resource for the researchers and graduate students in this area.

Let $k$ be a field of characteristic $\not=2$ and
$R_{n,m}=k[x_1,\ldots,x_n,y_1,\ldots,y_m]/(\sum_1^n x_i^2 -\sum_1^m
y_j^2-1)$. Then we will prove the following results.

\begin{theorem}
Assume that $x^2+y^2+z^2=0$ has only trivial zero in $k^3$
(equivalently the quaternion algebra $\frac {(-1,-1)}{k}$ is a
division algebra over $k$). Then $\wt K_0(R_{n,0})$ and $\wt
K_0(R_{0,m})$ are periodic of period $8$. More precisely,

$(1)$ $\wt
K_0(R_{n,0})$ is $\BZ, \BZ/2\BZ$ or $0$ depending on whether
$n$ is $\{1,5\}$ modulo $8$, $\{2,3\}$ modulo $8$ or $\{0,4,6,7\}$ modulo $8$.

$(2)$ $\wt K_0(R_{0,m})$ is $\BZ,\BZ/2\BZ$ or $0$ depending on whether
$m$ is $\{3,7\}$ modulo $8$, $\{5,6\}$ modulo $8$ or $\{0,1,2,4\}$
modulo $8$.

$(3)$ $\wt K_0(R_{n,m})=\wt K_0(R_{n-m,0})$ if $n\geq
m$ and $ \wt K_0(R_{n,m})=  \wt K_0(R_{0,m-n})$ if $n< m$.
\end{theorem}

\begin{theorem}
Assume $\sqrt {-1}\in k$. Then $\wt K_0(R_{n,m})$ is $\BZ$ or $0$
depending on whether $n+m$ is odd or even.
\end{theorem}

\begin{theorem}
Assume that $\sqrt {-1}\notin k$ and $-1$ is a sum of two squares in
$k$ (equivalently, the quaternion algebra $(\frac {-1,-1}{k})$ is not
a division algebra over $k$).  Then $\wt K_0(R_{0,n})$ and $\wt K_0(R_{n,0})$
are periodic of period $4$. More precisely,

$(i)$ $\wt K_0(R_{0,n})=\BZ,\BZ/2\BZ$ or $0$ depending on whether $n$ is
$3$ modulo $4$, $2$ modulo $4$ or $\{0,1\}$ modulo $4$ respectively.

$(ii)$ $\wt K_0(R_{n,0})=\BZ,\BZ/2\BZ$ or $0$ depending on whether $n$ is
$1$ modulo $4$, $0$ modulo $4$ or $\{2,3\}$ modulo $4$ respectively.

$(iii)$ $\wt K_0(R_{n,m})=\wt K_0(R_{n-m,0})$ if $n\geq m$ and $\wt
K_0(R_{n,m})=\wt K_0(R_{0,m-n})$ if $n< m$.
\end{theorem}


\section{Preliminaries}

We will recall some results from \cite{swan} for later use.  Let $k$
be a field of characteristic $\not =2$ and let $q(x_1,\ldots,x_n)$ be
a non-degenerate quadratic form over $k$.  Let
$R(q)=k[x_1,\ldots,x_n]/(q-1)$ and let $C(q)$ be the Clifford algebra
of $q$. If $q=\sum_1^n a_ix_i^2$, $a_i\in k$, then $C(q)$ is generated
by $e_1,\ldots,e_n$ with relations $e_ie_j+e_je_i=0$ for $i\not= j$
and $e_i^2=a_i$. The elements $e_{i_1}\ldots e_{i_r}$ with $1\leq
i_1<\ldots < i_r\leq n$ form a $k$-base for $C(q)$.

If $q=a_1x_1^2+\ldots+a_nx_n^2$, then $det\, q =a_1\ldots a_n$ and $ds
\,q=(-1)^{n(n-1)/2} det\, q$.  A binary quadratic form is called
hyperbolic if it has the form $h(x,y)=x^2-y^2$. By a linear change of
variables this is equivalent to $h'(x,y)=xy$.

\begin{lemma}\label{8.1}
(\cite{swan}, 8.1 and 8.2) If $b$ is a binary quadratic form, then
$C(b\perp q)\iso C(b)\ot C((ds\,b)q)$. In particular, if $h$ is
hyperbolic, then $C(q\perp h)\iso C(q)\ot C(h)$.
\end{lemma}

\begin{lemma}\label{8.3}
(\cite{swan}, 8.3) $(a)$ If $q$ has even rank, then $C(q)$ is central
simple over $k$ and is a tensor product of quaternion algebras.

$(b)$ If $q$ has odd rank, then $(i)$ if $\sqrt {ds\, q}\in k$, then
$C(q)=A\times A$, where $A$ is central simple over $k$ and is a tensor
product of quaternion algebras, $(ii)$ otherwise $C(q)$ is simple with
center $k(\sqrt {ds\, q})$ and is a tensor product of its center with
quaternion algebras over $k$.
\end{lemma}

It follows from (\ref{8.3}) that all simple $C(q)$-modules have the
same dimension over $k$. Denote this dimension by $d(q)$.

\begin{lemma}
(\cite{swan}, Lemma 8.4)
$(a)$ $d(q\perp 1)$ is either $d(q)$ or $2d(q)$.

$(b)$ If $C(q)=A\times A$, then $d(q\perp 1)=2d(q)$.
\end{lemma}

See \cite{swan} for the definition of $ABS(q)$.

\begin{proposition}\label{1}
(\cite{swan}, Proposition 8.5) $(a)$ If $C(q)=A\times A$, i.e rank of
$q$ is odd and $\sqrt {ds \,q}\in k$, then $ABS(q)=\BZ$ generated by
either of the simple $C(q)$-modules.

$(b)$ If $C(q)$ is simple, then $(i)$ $ABS(q)=0$ if $d(q\perp 1)=d(q)$
and $(ii)$ $ABS(q)=\BZ/2\BZ$ if $d(q\perp 1) =2d(q)$.
\end{proposition}

We state the following result of Swan (\cite{swan1}, Corollary 10.8)

\begin{theorem}\label{2}
Assume that $R$ is regular, $1/2\in R$ and $q\perp <-1>$ is a
non-singular quadratic form. Then $ABS(q) \iso K_0(R(q))/K_0(R)$,
where $R(q)=R[x_1,\ldots,x_n]/(q-1)$.

In particular, if $R$ is a field, then $ABS(q)\iso \wt K_0(R(q))$.
\end{theorem}

Using (\ref{1} and \ref{2}), we get the following result which will be
used later.

\begin{theorem}\label{12}
Let $k$ be a field of characteristic $\not= 2$ and let
$q(x_1,\ldots,x_n) \perp <-1>$ be a non-singular quadratic form over
$k$. Write $R(q)=k[x_1,\ldots,x_n]/(q(x_1,\ldots,x_n)-1)$. Then we
have

$(i)$ If $C(q)=A\times A$ (i.e. rank of $q$ is odd and $\sqrt {ds\, q}
\in k$), then $\wt K_0(R(q))=\BZ$.

$(ii)$ If $C(q)$ is simple, then $(a)$ $\wt K_0(R(q))=0$ if $d(q\perp
1)=d(q)$ and $(b)$ $\wt K_0(R(q))=\BZ/2\BZ$ if $d(q\perp 1) =2d(q)$.
\end{theorem}


\section{Main Theorem}

\subsection{$-1$ is not a sum of two squares in $k$}

Let us recall the following well known result (see \cite{P}, p. 15).
Let $k$ be a field and let $a,b\in k$. Then the quaternion algebra
$\frac{(a,b)}{k}$, which is a $k$-algebra defined by $i$ and $j$ with
relations $i^2=a$, $j^2=b$ and $ij+ji=0$, is a division algebra if and
only if $x^2=ay^2+bz^2$ has only trivial zero.

In this section we will assume that $k$ is a field such that
 $x^2+y^2+z^2=0$ has only trivial zero in $k^3$ which is same as the
 quaternion algebra $\frac{(-1,-1)}{k}$ is a division algebra over $k$
 (e.g. any real field). We denote the division algebra
 $\frac{(-1,-1)}{k}$ by $\CH$. Let $\CC$ be the subalgebra of $\CH$
 generated by $i$ over $k$. Then $\CC=k[x]/(x^2+1)$ is a field.

The following is a well known result. We will give proof for
completeness. Recall that $F(n)$ denote the algebra of $n\times n$
matrices over $F$.

\begin{lemma}
If $F$ denote one of $k$, $\CC$ or $\CH$, then we have the following
identities $(i)$ $F(n)\iso k(n)\ot_k F$, $(ii)$ $k(n)\ot_k k(m)\iso
k(nm)$, $(iii)$ $\CC\ot_k \CC\iso \CC\op \CC$, $(iv)$ $\CH\ot_k
\CC\iso \CC(2)$, $(v)$ $\CH\ot_k \CH \iso k(4)$.

In particular, when $k=\BR$ the field of real numbers, then $\CC=\BC$
and $\CH=\BH$.
\end{lemma}

\begin{proof}
$(i)$ and $(ii)$ are straightforward.

$(iii)$ The map $\CC\op \CC \ra \CC\ot_k \CC$ defined by $(1,0)
\mapsto 1/2(1\ot 1+i\ot i)$ and $(0,1)\mapsto 1/2(1\ot 1- i\ot i)$ is
an isomorphism.

$(iv)$ Since $\CH$ is a $\CC$ vector space under left multiplication,
the map $\pi : \CC\times \CH \ra \Hom_{\CC}(\CH,\CH)$ defined by
$\pi_{y,z}(x)=yx\ol z$ is $k$-bilinear, where $y\in \CC$, $x,z\in \CH$
and $\ol z=a1-bi-cj-dij$ is the conjugate of $z=a1+bi+cj+dij$ with
$a,b,c,d\in k$. Hence, we get a $k$-linear map $\pi : \CC\ot_k \CH \ra
\Hom_{\CC}(\CH,\CH)$. Since $\pi_{y,z}\circ \pi_{y^\prime,
z^\prime}=\pi_{yy^\prime,zz^\prime}$, the map $\pi$ is an $k$-algebra
homomorphism. Further, it is easy to see that $\pi$ is
injective. Since $\Hom_{\CC}(\CH,\CH)\iso \CC(2)$, we get $\dim_k
\CC\ot_k \CH =8 =\dim_k \CC(2)$ (note that $\dim_{\CC}
\CC(2)=4$). Hence $\pi$ is an isomorphism.

$(v)$ Define a map $\pi: \CH\times \CH \ra \Hom_{k}(\CH,\CH)$ by
$\pi_{y,z}(x)=yx\ol z$, where $y,x,z\in \CH$. Then $\pi$ is
$k$-bilinear. Hence it induces a $k$-linear map $\pi : \CH\ot_k \CH
\ra \Hom_k(\CH,\CH)$, which is an algebra homomorphism
($\pi_{y,z}\circ
\pi_{y^\prime,z^\prime}=\pi_{yy^\prime,zz^\prime}$). Further, $\pi$ is
injective. Since both sides are vector spaces of dimension $16$ over
$k$, $\pi$ is an isomorphism. Note that $\Hom_k(\CH,\CH)\iso
k(4)$. This proves the result.  $\hfill \square$
\end{proof}

Let $q_n=-\sum_1^n x_i^2$ and $q_n'=\sum_1^n x_i^2$ be quadratic forms
over $k$.  We write $C_n$ and $C_n'$ for the Clifford algebras
$C(q_n)$ and $C(q_n')$ respectively. Then we have the following
result. In \cite{ABS}, it is proved for $k=\BR$, but the same proof
works over any field $k$.

\begin{proposition}\label{4.2}
(\cite{ABS}, Proposition 4.2) There exist isomorphisms $C_n \ot_{k}
C_2' \iso C_{n+2}'$ and $C_n'\ot_{k} C_2\iso C_{n+2}$.
\end{proposition}

It is easy to see that $C_1=\CC$, $C_2=\CH$, $C_1^\prime =k\op k$ and
$C_2^\prime =k(2)$.  Using (\ref{4.2}), we get that

$$\begin{array}{|c|c|c|c|c|}
\hline

n & C_n & C_n' & d(q_n) & d(q_n')\\ \hline

1  &  \CC   & k\op k      &  2       &1\\
2  &  \CH   & k(2)      &     4    &2\\
3  &  \CH\op \CH    & \CC(2)      &4         &4\\
4  &  \CH(2)    & \CH(2)      &     8    &8\\
5  &  \CC(4)    &  \CH(2)\op \CH(2)     &8         &8\\
6  &  k(8)    &  \CH(4)     &       8  &16\\
7  &  k(8) \op k(8)   & \CC(8)      & 8        &16\\
8  &  k(16)   &  k(16)     &       16  &16 \\ \hline
\end{array}$$

Note that $C_4\iso C_4'$, $C_{n+4}\iso C_n \ot_{k} C_4$, $C_{n+8}\iso
C_n\ot C_8$. Further $C_8\iso k(16)$. Hence if $C_n=F(m)$, then
$C_{n+8}\iso F(16 m)$. Similarly, if $C_n^\prime =F(m)$, then
$C_{n+8}^\prime=F(16 m)$.

Let $h=x^2-y^2$ be the hyperbolic quadratic form over $k$. Then the
Clifford algebra $C(h)\iso k(2)$. From (\ref{8.1}), the Clifford
algebra of $h^r=h\perp \ldots \perp h$ ($r$ times) is $C(h^r)=k(2)\ot
\ldots \ot k(2) \iso k(2^r)$.  Now, if $q=q_1 \perp h^r$ and
$C(q_1)=F(m)$, then $C(q)=F(m)\ot k(2^r)\iso F(2^rm)$.

Since $q_k \perp 1 = q_{k-1} \perp h$, the Clifford algebra $C(q_k
\perp 1)=C(q_{k-1}) \ot k(2)$.  Write $q_k \perp 1$ as $\wt q_k$.
Further $q_n^\prime \perp 1 = q^\prime_{n+1}$. Hence the Clifford algebra
$C(q^\prime_n \perp 1)=C(q^\prime_{n+1})$ and $d(q_k^\prime \perp
1)=d(q_{k+1}^\prime)$. 

Write $s=16^r$. Then we have the following:

$$\begin{array}{|c|c|c|c|c|c|c|}
\hline

n & C_{8r+n} & C_{8r+n}' & C(\wt q_{8r+n}) & d(q_{8r+n}) & d(q_{8r+n}') & 
d(\wt q_{8r+n})\\ \hline

1  &  \CC(s)   & k(s)^2 &  k(2s)   &  2s      &s    &2s \\
2  &  \CH(s)   & k(2s)   & \CC(2s)   &     4s   &2s      &4s \\
3  &  \CH(s)^2    & \CC(2s)  & \CH(2s)    &4s   &4s &8s \\
4  &  \CH(2s)    & \CH(2s)   & \CH(2s)^2   &  8s  &8s &8s \\
5  &  \CC(4s)    &  \CH(2s)^2 & \CH(4s)  &8s  &8s &16s\\
6  &  k(8s)    &  \CH(4s)   & \CC(8s)  &       8s  &16s &16s \\
7  &  k(8s)^2 & \CC(8s)  &k(16s)    & 8s       &16s &16s \\
8  &  k(16s)   &  k(16s)     & k(16s)^2     & 16s  &16s &16s \\ \hline
\end{array}$$

Using (\ref{12}), we get the following result.

\begin{theorem}
Let $k$ be a field of characteristic $\not=2$ such that
$x^2+y^2+z^2=0$ has only trivial zero in $k^3$.  Let
$R_{0,n}=R(q_n)=k[x_1,\ldots,x_n]/(-\sum_1^n x_i^2-1)$ and
$R_{n,0}=R(q_n')=k[x_1,\ldots,x_n]/(\sum_1^n x_i^2-1)$. 
Then we have the following:

$(1)$ $\wt
K_0(R_{n,0})$ is $\BZ, \BZ/2\BZ$ or $0$ depending on whether
$n$ is $\{1,5\}$ modulo $8$, $\{2,3\}$ modulo $8$ or $\{0,4,6,7\}$ modulo $8$.

$(2)$ $\wt K_0(R_{0,m})$ is $\BZ,\BZ/2\BZ$ or $0$ depending on whether
$m$ is $\{3,7\}$ modulo $8$, $\{5,6\}$ modulo $8$ or $\{0,1,2,4\}$
modulo $8$.

\end{theorem} 

Now, let $n,m$ be positive integers and consider the quadratic form
$Q_{n,m}(x_1,\ldots,x_n,y_1,\ldots,y_m)=\sum_1^n x_i^2- \sum_1^m
y_i^2$. First assume that $n\geq m$. Then $Q_{n,m} \iso q_{n-m}^\prime
\perp h^m$ and the Clifford algebra $C(Q_{n,m})$ is isomorphic to
$C_{n-m}^\prime \ot k(2^m)$. Hence $d(Q_{n,m})=d(q_{n-m}^\prime) 2^m$.
Further, $Q_{n,m} \perp 1 \iso q_{n-m+1}^\prime \perp h^m$ and
$d(Q_{n,m}\perp 1)=d(q_{n-m+1}^\prime) 2^m$.  Hence $d(Q_{n,m}\perp
1)/d(Q_{n,m}) = d(q_{n-m}^\prime)/d(q_{n-m+1}^\prime)$. 

Now assume that $n <m$. Then $Q_{n,m} \iso q_{m-n} \perp h^n$ and
$Q_{n,m}\perp 1 \iso q_{m-n-1} \perp h^{n+1}$. Further the Clifford
algebra $C(Q_{n,m})\iso C(q_{m-n})\ot k(2^n)$ and $C(Q_{n,m}\perp 1)
=C(q_{m-n-1})\ot k(2^{n+1})$. Hence, $d(Q_{n,m})=d(q_{m-n}) 2^n$ and
$d(Q_{n,m}\perp 1)=d(q_{m-n-1}) 2^{n+1}$.  The quotient
$d(Q_{n,m}\perp 1)/d(Q_{n,m})$ is equal to $2
d(q_{m-n-1})/d(q_{m-n})$.

By (\ref{12}), we get the following:

\begin{theorem}\label{t1}
Let $k$ be a field of characteristic $\not=2$ such that
$x^2+y^2+z^2=0$ has only trivial zero in $k^3$. Let
$R(Q_{n,m})=k[x_1,\ldots,x_n,y_1,\ldots,y_m]/(Q_{n,m}-1)$. Then $\wt
K_0(R(Q_{n,m}))$ is same as $\wt K_0(R(q_{n-m}^\prime))$ when $n\geq
m$ and $\wt K_0(R(q_{m-n}))$ when $n <m$.
\end{theorem}

\begin{remark}
We note that the following classical result generalizes (\ref{t1})
(see \cite{swan}, 10.1). Let $k$ be a field of characteristic $\not=2$
and let $f\not=0 \in k[x_1,\ldots,x_n]$. Let $A=k[x_1,\ldots,x_n]/(f)$
and $B=k[x_1,\ldots,x_n,u,v]/(f+uv)$. Then $\wt G_0(A) \iso \wt
G_0(B)$ and for a regular ring $R$, we know that $\wt G_0(R) \iso \wt
K_0(R)$. We give direct proof and compute $\wt K_0(R(q_n))$.
\end{remark}


\subsection{$\sqrt{-1}\in k$}

Let $k$ be a field of characteristic $\not=2$ such that $\sqrt {-1}
\in k$.  Let $q_n(x_1,\ldots,x_n)=-(x_1^2+\ldots+x_n^2)$ and
$q_n^\prime(x_1,\ldots,x_n)=x_1^2+\ldots+x_n^2$ be quadratic forms
over $k$. Let $C_n$ and $C_n^\prime$ be the Clifford algebras of $q_n$
and $q_n^\prime$ over $k$ respectively. Then $C_n \iso C_n^\prime$.
Further, using (\ref{4.2}), we get $C_{n+2} \iso C_n\ot C_2$. Since
$C_1 =k\op k$ and $C_2=k(2)$, we get $C_{2n}=k(2^n)$ and
$C_{2n+1}=k(2^n)\op k(2^n)$. Therefore, by (\ref{12}), we get
the following result.

\begin{theorem}
Let $k$ be a field of characteristic $\not=2$ such that $\sqrt {-1}\in
k$. Let $q_n=x_1^2+\ldots+x_n^2$ and
$R(q_n)=k[x_1,\ldots,x_n]/(q_n-1)$. Then $\wt K_0(R(q_{2n}))=0$ and
$\wt K_0(R(q_{2n+1}))=\BZ$.
\end{theorem}


\subsection{$-1$ is a sum of two squares and $\sqrt {-1}\notin k$}

Let $k$ be a field of characteristic $\ne 2$ such that $\sqrt
{-1}\notin k$ but $-1$ is a sum of two squares in $k$
(i.e. $x^2+y^2+z^2=0$ has a non-trivial zero in $k^3$).
Let $q_n(x_1,\ldots,x_n)=-(x_1^2+\ldots+x_n^2)$ and
$q_n^\prime(x_1,\ldots,x_n)=x_1^2+\ldots+x_n^2$ be quadratic forms
over $k$. Let $C_n$ and $C_n^\prime$ be the Clifford algebras of $q_n$
and $q_n^\prime$ over $k$ respectively.

We denote the field $k[x]/(x^2+1)$ by $\CC$. We recall the following
well known result: A quaternion algebra $(\frac {a,b} k)$ is isomorphic
to $M_2(k)$ if and only if it is not a division algebra. Then it is
easy to see that $C_1=\CC$, $C_1'=k\op k$, $C_2=k(2)=C_2'$. Further,
$C_3=C_1'\ot C_2=k(2)\op k(2)$, $C_3'=C_1\ot C_2'= \CC(2)$ and
$C_4=C_2'\ot C_2=k(4)=C_4'$.

For $n=4r+i$, where $i\in \{1,2,3,4\}$, we have $C_n=C_{n-2}'\ot C_2 =
C_{n-4}\ot C_2'\ot C_2=C_{n-4}\ot C_4=C_{n-4}\ot k(4)=\ldots =C_i\ot
k(4^r)$. Similarly, $C'_n=C_i'\ot k(4^r)$. 

Write $s=4^r$. Then we have the followings:
$$\begin{array}{|c|c|c|c|c|c|c|}
\hline

n & C_{4r+n} & C_{4r+n}' & C(q_{4r+n}\perp 1) & d(q_{4r+n}) & d(q_{4r+n}') & 
d(q_{4r+n}\perp 1)\\ \hline

1  &  \CC(s)   & k(s)^2 &  k(2s)   &  2s      &s    &2s \\
2  &  k(2s)   & k(2s)   & \CC(2s)   &     2s   &2s      &4s \\
3  &  k(2s)^2    & \CC(2s)  & k(4s)    &2s   &4s &4s \\
4  &  k(4s)    & k(4s)   & k(4s)^2   &  4s  &4s &4s\\
\hline
\end{array}$$

By (\ref{12}), we get the following result.

\begin{theorem}
Let $k$ be a field of characteristic $\not=2$ such that $\sqrt
{-1}\notin k$ and $-1$ is a sum of two squares in $k$.  Let
$R_{n,m}=k[x_1,\ldots,x_n,y_1,\ldots, y_m]/(\sum_1^n x_i^2 -\sum_1^m
y_j^2-1)$. Then

$(i)$ $\wt K_0(R_{0,n})=\BZ,\BZ/2\BZ$ or $0$ depending on whether $n$ is
$3$ modulo $4$, $2$ modulo $4$ or $\{0,1\}$ modulo $4$ respectively.

$(ii)$ $\wt K_0(R_{n,0})=\BZ,\BZ/2\BZ$ or $0$ depending on whether $n$ is
$1$ modulo $4$, $0$ modulo $4$ or $\{2,3\}$ modulo $4$ respectively.

$(iii)$ $\wt K_0(R_{n,m})=\wt K_0(R_{n-m,0})$ if $n\geq m$ and $\wt
K_0(R_{n,m})=\wt K_0(R_{0,m-n})$ if $n< m$.
\end{theorem}


\section{Some Auxiliary Results}

Let $A=\BR[x_0,\ldots,x_n]/(a_0x_0^2+\ldots+a_nx_n^2-b)$ with
$a_i,b\in \BR$ and let $E(A)$ be the Euler class group of $A$ with
respect to $A$ (see \cite{BR1} for definition). Let $E^{\BC}(A)$ be
the subgroup of $E(A)$ generated by all the complex maximal ideals of
$A$.  Then $E^{\BC}(A)=0$, since by (\cite{MR}, Lemma 4.2), all the
complex maximal ideals of $A$ are generated by $n$ elements. Hence by
(\cite{MS}, Theorem 2.3), we get the following results: $(i)$
$E(A)\iso E(\BR(X))$, where $X=\Spec (A)$ and $\BR(X)$ is the
localization $A_S$ of $A$ with $S$ as the set of all elements of $A$
which do not have any real zero and $(ii)$ $CH_0(A)\iso
CH_0(\BR(X))$. Further, there is a natural surjection from $E(A)$ to
$CH_0(A)$.

$(1)$ Assume that $A=\BR[x_0,\ldots,x_n]/(x_0^2+\ldots+x_n^2+1)$. Then
$A$ has no real maximal ideal and hence $E(A)=E^{\BC}(A)=0$ and hence
$CH_0(A)=0$. For $A=\BR[x_0,\ldots,x_n]/(x_0^2+\ldots+x_n^2-1)$, it is
known that $E(A)=\BZ$ and $CH_0(A)=\BZ/2$.

$(2)$ Assume $A=\BR[x_0,\ldots,x_n]/(\sum_0^m x_i^2 -\sum_{m+1}^n
x_i^2 -1)$ with $m<n$ and $X=\Spec (A)$. Then $X(\BR)$ has no compact
connected component. Hence by (\cite{BDM}, Theorem 4.21),
$E(\BR(X))=0$. From above, we get $E(A)=0$ and $CH_0(A)=0$.

$(3)$ In general, let
$A=\BR[x,y,z_1,\ldots,z_n]/(xy+f(z_1,\ldots,z_n))$ and let $X=\Spec
(A)$. Then $X(\BR)$ has no compact connected component. All the
connected components of $X(\BR)$ is unbounded. For this, note that if
$(a,b,c_1,\ldots,c_n)\in X(\BR)$, then $f(c_1,\ldots,c_n)=-ab$ and if
$(x_0,y_0)$ is any point on the hyperbola $xy=ab$, then
$(x_0,y_0,c_1,\ldots,c_n)\in X(\BR)$. Hence, by (\cite{BDM}, Theorem
4.21), the Euler class group of $\BR(X)$, namely $E(\BR(X))=0$ and we
get that $E(A)=E^{\BC}(A)$, the group generated by all the complex
maximal ideals of $A$. Using (\cite{MS}, Theorem 2.3), we get that
$E(A)\iso CH_0(A)$. Further, it is known (see \cite{BR1}, Theorem 5.5)
that for a smooth affine domain $A$ of dimension $\geq 2$ over $\BR$,
$CH_0(A) \iso E_0(A)$, the weal Euler class group of $A$. Hence
$E(A)\iso E_0(A)\iso CH_0(A)$ and $E(A)$ is generated by complex
maximal ideals of $A$. In particular, if all the complex maximal
ideals of $A$ are generated by $n$ elements, then $E(A)=0$ as is the
case in $(2)$ above.

{}

\end{document}